\documentclass[a4paper]{article}
\usepackage{color}
\usepackage{caption}
\usepackage{float}
\usepackage{amsmath} 
\usepackage{calc}
\newlength{\depthofsumsign}
\usepackage{graphicx,epstopdf}
\usepackage{amsfonts,amsmath,amssymb,amsthm}
\usepackage{graphicx}
\usepackage{mathrsfs}
\usepackage[title, titletoc, header]{appendix}

\newtheorem{definition}{Definition}[section]

\newtheorem{theorem}{Theorem}[section]

\newtheorem{example}{Example}[section]
\usepackage{amsmath}

\usepackage{graphicx,epstopdf}

\newtheorem{remark}{Remark}[section]

\begin{document}

\begin{center}
\large{ \textbf{Positive constrained approximation via RBF-based partition of unity method}}
\end{center}

\begin{center}
Alessandra De Rossi and Emma Perracchione 
\end{center}

\begin{center}
Department of Mathematics \lq\lq G. Peano\rq\rq, University of Torino, via Carlo Alberto 10, I--10123 Torino, Italy
\end{center}
\vskip 0.5cm

\textbf{Abstract.}
In this paper, we discuss the problem of constructing Radial Basis Function (RBF)-based Partition of Unity (PU) interpolants that are positive if data values are positive. 
More specifically, we compute positive local approximants by adding up several constraints to the interpolation conditions. This approach, considering a global approximation problem and Compactly Supported RBFs (CSRBFs), has been previously proposed in \cite{WU10a}. Here, the use of the PU technique enables us to intervene only locally and as a consequence to reach a better accuracy. This is also due to the fact that we select the \emph{optimal} number of positive constraints by means of an \emph{a priori} error estimate and we do not restrict to the use of CSRBFs.  Numerical experiments and applications to population dynamics are provided to illustrate the effectiveness of the method in applied sciences.

\section{Introduction}
\label{Intro}

Given a set of samples, the scattered data interpolation problem consists in finding an approximating function that matches the given measurements at their corresponding locations. Furthermore, dealing with applications, we often have additional properties, such as the non-negativity of the measurements, which we wish to be preserved during the interpolation process. Note that, since such property is known as \emph{positivity-preserving property}, to keep a common notation with existing literature, we use the term positive (instead of non-negative) function values or interpolants.

To preserve such property, mostly considering rational spline functions with $C^1$ or $C^2$ continuity, the recent research studies techniques which force the approximants to be positive. As example, the conditions under which the positivity of a cubic piece may be lost are investigated in \cite{Schmidt}. Moreover, in order to preserve the positivity, the use of bicubic splines, coupled with a technique based on adding extra knots, has been investigated in \cite{Asim1,Butt}. A positive fit is instead obtained by means of rational cubic splines in \cite{Hussain06,Hussain08}. The same authors also developed a positive surface construction scheme for positive scattered data arranged over triangular grids  \cite{Hussain10}.

Note that all the above mentioned methods depend on a mesh. However, the positivity-preserving problem is also well-known in the field of meshfree or meshless methods. They include Shepard-type approximants \cite{Renka,Shepard} and RBF interpolants \cite{Buhmann03,Wendland05}. While the  positivity-preserving problem has been widely investigated for the Modified Quadratic Shepard's (MQS) approximant \cite{Asim,Brodlie}, it remains a challenging computational issue for RBF interpolation. Indeed, even if such meshfree approach has been extensively studied in the recent years, especially focusing on the stability of the interpolant  \cite{Demarchi15,Fornberg11}, not a lot of effort has been addressed to construct positivity-preserving approximants. Such problem has been studied only in particular and well-known cases; for instance, in \cite{Utreras} it is analyzed for the thin plate spline. Other methods, which follow from the quasi-interpolation formula given in \cite{WU94}, have been proposed and effectively performed to construct positive approximants of positive samples  \cite{Wang}.

Focusing on RBFs, our scope consists in preserving the positivity of the PU interpolant for a wider family of kernels. In \cite{WU10a} a global positive RBF approximant is constructed by adding up to the interpolation conditions several positive constraints and considering CSRBFs. Even if the optimal number of constraints is not investigated, the results are promising and show that this technique has a better accuracy than the Constrained MQS (CMQS) approximant. However, since a global interpolant is used, adding up other constraints to preserve the positivity implies that the shape of the curve/surface is consequently globally modified. As pointed out in \cite{WU10a}, this might lead to a considerable decrease of the quality of the approximating function in comparison with the unconstrained CSRBF interpolation.

Thus here, in order to avoid such drawback, focusing on 2D data sets, the PU method is performed by imposing positive constraints on the local RBF interpolants. Such approach enables us to consider constrained interpolation problems only in those PU subdomains which do not preserve the required property. This leads to an accurate method compared with existing techniques \cite{Asim,Brodlie,WU10a}. Moreover, in contrast with \cite{WU10a}, we can consider truly large data sets.

Specifically, in order to construct the Positive Constrained PU (PC-PU) approximant, following \cite{WU10a}, we locally impose several positive constraints and  we reduce to solve an optimization problem. The number of constraints is properly selected by means of an a priori error estimate.   This is a fundamental step to maintain a good accuracy of the fit. Moreover, differently from \cite{WU10a}, a wider family of RBFs (not only compactly supported) is considered. The main disadvantage of using CSRBFs is that they introduce large errors in the area in which the interpolant is negative. This is due to the fact that the shape of the fit is modified only within the support of the CSRBF and thus neighbouring points are not taken into account. Numerical evidence shows that the use of infinitely smooth globally defined RBFs leads to an improvement in this direction.
Moreover, in order to stress the effectiveness of the proposed technique in applied sciences, we investigate several applications to biomathematics. Finally, comparisons with the unconstrained PU interpolant, with the Shepard's method and with the one outlined in \cite{WU10a} are carried out.

The guidelines of the paper are as follows. In Section \ref{positive_interpol} we investigate the positivity  of the PU approximant by considering extra positive constraints. The computational aspects  of such method are analyzed in Section \ref{locations}. Numerical experiments and applications to population dynamics are shown in Section \ref{NE_POS} and \ref{app_erbivori}, respectively. Finally, Section \ref{cr} deals with conclusions and future work.

\section{Positive approximation of positive data values}
\label{positive_interpol}

In addition to the assigned interpolation conditions, our goal consists in considering several positive constraints. This allows to preserve the positivity of the measurements. Therefore, our approach turns out to be meaningful especially in applications, indeed in order to avoid the violation of  biological or  physical  measurements, a positive fit is often necessary. Thus, before discussing the proposed method and without loosing generality, we consider an example of univariate interpolation to illustrate the scope of such numerical tool and motivate the reader. 

\begin{example}[Motivations and targets]
	One of the most common tumor, affecting mainly men over sixty years old, is the prostate cancer. Luckily, this disease has a very slow growth and a reliable biomarker after prostatectomy, the so-called Prostate Specific Antigen (PSA), suitable for the early diagnosis. In particular, if its value is larger than $0.2$ ng/mL,  a relapse (a local or distal metastasis) occurs.
	
	Clinical data of prostatectomized patients are available in \cite{Gabriele}. They are used to investigate the evolution of the relapse via mathematical models \cite{Stura2Pop}. 	Therefore, in order to validate such models, a data fitting results  essential \cite{Perracchione}. For such scope, one can reconstruct the PSA curve with the standard PU interpolant. Anyway, the positivity of the PSA values is not always preserved through the interpolation process, violating the biological constraint. In order to avoid this problem, we propose a novel technique, namely the PC-PU method, which allows to preserve the positivity property, see Figure \ref{psa}.
	
	\begin{figure}[ht!]
		\begin{center}
			\makebox[\textwidth]{
				\begin{minipage}[b]{0.5\textwidth}
					\centering
					\includegraphics[width=\textwidth]{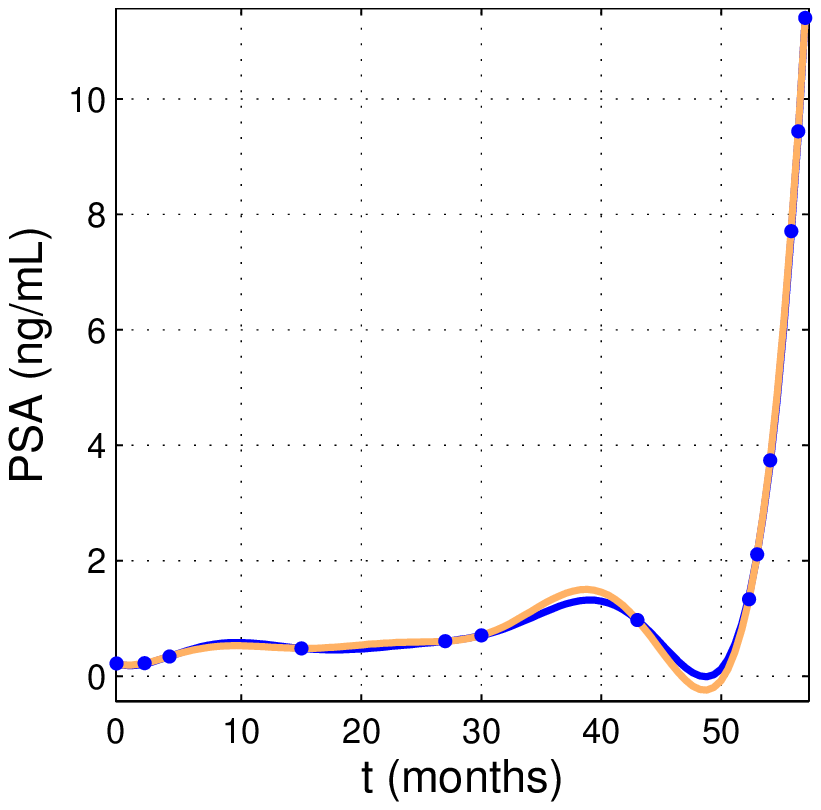} 
				\end{minipage}
				\begin{minipage}[b]{0.5\textwidth}
					\centering
					\includegraphics[width=\textwidth]{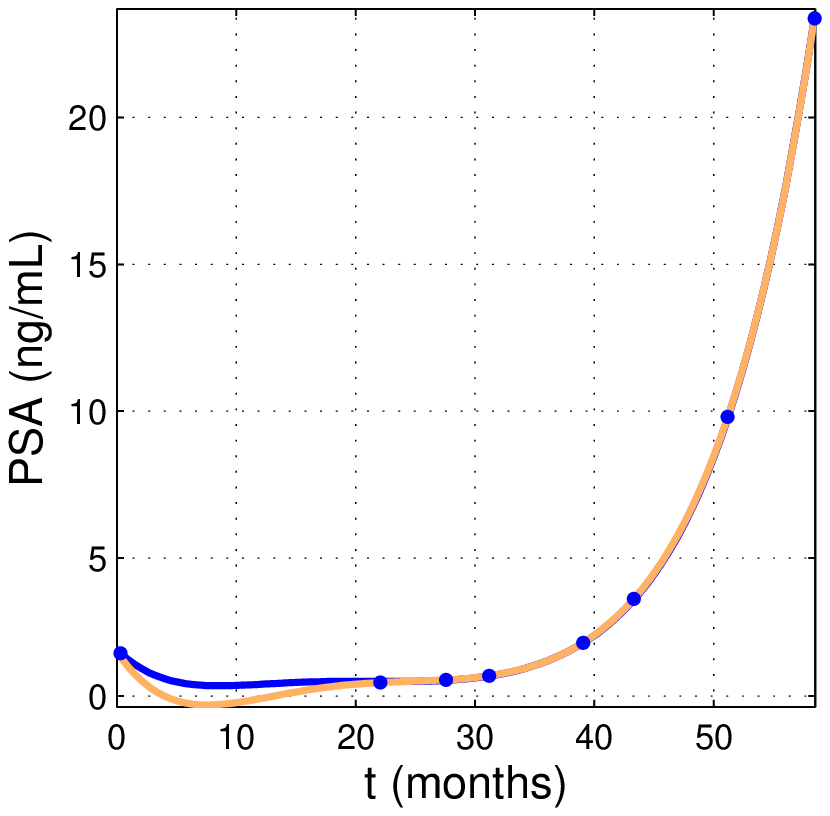} 
				\end{minipage}  
			}				
		\end{center} 
		\caption{Examples of curves fitting PSA values (plotted with blue dots). The classical PU interpolant is plotted in orange and the  PC-PU approximant in blue.}
		\label{psa}
	\end{figure}
\end{example}

\subsection{The partition of unity method}
Let us consider a bounded set $\Omega \subseteq \mathbb{R}^2$.
Given a set  of $N$ distinct data points (or data sites or nodes) ${ \cal X}_N = \{(x_i,y_i), i=1, \ldots, N \} \subseteq \Omega$ and a set of data values (or measurements  or function values) ${\cal F}_N = \{f_i=f(x_i,y_i), i = 1, \ldots ,N \}$, the  PU method decomposes the domain $\Omega$ into $d$ subdomains $ \Omega_j$, $j=1, \ldots, d$, such that $\Omega  \subseteq \bigcup_{j=1}^{d} \Omega_j$ \cite{Babuska97,Fasshauer,Melenk96,Wendland02a}. In literature, the subdomains $\Omega_j$ are  supposed to be circular patches of the same radius $\delta$  covering the domain $\Omega$. 

According to \cite{Wendland02a}, associated with these subdomains, we choose a $p$-stable partition of unity, i.e a family of nonnegative functions $ \{ W_j \}_{j=1}^{d}$,  with $W_j \in C^p ( \mathbb{R}^2)$, such that
\begin{itemize}
	\item[i.] supp$ ( W_j ) \subseteq  \Omega_j $,
	\item[ii.] $ \sum_{j=1}^{d} W_j( x,y )=1$ on $ \Omega$,
	\item[iii.] $||D^{ \beta} W_j ||_{L^{ \infty} ( \Omega_j)} \leq \frac{C_{ \beta} }{ \xi_j^{ | \beta|}},$ $  \forall \beta \in \mathbb{N}^{2}: | \beta |	 \leq p,$
	where $C_{ \beta} > 0$ is a constant and $\xi_j= \textrm{diam}(\Omega_j)$.
\end{itemize}

Then, the global interpolant assumes the form
\begin{equation} 
{\cal I}(x,y)= \sum_{j=1}^{d} R_j( x,y) W_j (x,y), \quad (x,y) \in \Omega,
\label{intg}
\end{equation}
where $R_j$ defines a local interpolant on each PU subdomain $\Omega_j$ and $W_j: \Omega_j \longrightarrow \mathbb{R}$. 
\begin{remark}
	Since the functions $W_j$, $j=1, \ldots, d$, form a partition of unity, if the local fits $R_j$, $j=1,\ldots,d$, satisfy the interpolation conditions, then the global PU approximant inherits the interpolation property, i.e.
	\begin{align*}
	{\cal I}(x_i,y_i ) = \sum_{j=1}^{d} R_j(x_i ,y_i) W_j ( x_i ,y_i) = \sum_{j \in I(x_i,y_i)}  f_i W_j ( x_i ,y_i) = f_i,
	\end{align*}
	where
	$I(x_i,y_i)= \{ j / (x_i,y_i) \in \Omega_j \}.$
\end{remark}
In what follows, the local interpolants $R_j$ are defined as  linear combinations of RBFs centred at $(x_k,y_k)$,  with  $(x_k,y_k) \in { \cal X}_{N_j}={ \cal X}_N \cap \Omega_j$, for $k=1, \ldots, N_j$, i.e.
\begin{equation}
R_j (x,y)= \sum_{k=1}^{N_j} c_k \phi_{\varepsilon}^{(k)} (x,y), \quad (x,y) \in \Omega_j,
\label{rad1}
\end{equation}
where $N_j$ is the number of nodes on $\Omega_j$ and $\phi_{\varepsilon}^{(k)} (x,y)$ is a RBF centred at $(x_k, y_k)$. A RBF depends on the so-called shape parameter $\varepsilon$ which governs the flatness of the function.

Among a large variety of  RBFs, we can differentiate between compactly supported and globally defined RBFs. As examples of these two major classes, we consider the Wendland's $C^2$ function and the Inverse MultiQuadric (IMQ). The latter is strictly positive definite and given by
\begin{equation}
\phi_{\varepsilon}^{(k)}  (x,y) = \dfrac{1}{\sqrt{1+ \varepsilon^2 r_k^2(x,y)}},
\label{IMQ}
\end{equation}
where $r_k^2(x,y)= (x -x_k)^2 + (y -y_k)^2$ is the square of the Euclidean distance from the centre $(x_k, y_k)$. Moreover, in case of ill-conditioning of the interpolation matrix, it might be advantageous to work with locally supported functions. Indeed, they lead to sparse linear systems. Wendland found a class of RBFs which are smooth, locally supported and strictly positive definite. For example, the  Wendland's $C^2$ function is defined as
\begin{equation}
\phi_{\varepsilon}^{(k)}  (x,y) =  (1- \varepsilon r_k(x,y))_{+}^4 (4\varepsilon r_k(x,y) +1),
\label{WEN}
\end{equation}
where  $( \cdot )_{+}$ denotes the truncated power function. Note that the shape parameter for a CSRBF  identifies the support of the function.

The coefficients $\{ c_k \}_{k=1}^{N_j}$ in \eqref{rad1} are determined by enforcing the $N_j$ local interpolation conditions
\begin{equation*}
R_j (x_i,y_i)= f_i, \quad i=1, \ldots, N_j,
\label{ls}
\end{equation*}
where $(x_i,y_i) \in { \cal X}_{N_j}$, and $f_i$ are the associated function values. Thus, the problem of finding the PU interpolant \eqref{intg} requires to solve $d$ linear systems of the form
\begin{equation*}
A^{(j)} \boldsymbol{c}= \boldsymbol{f} ,
\label{sys1}
\end{equation*}
where $\boldsymbol{c}= \left( c_1, \ldots, c_{N_j} \right)^T$, $  \boldsymbol{f} = \left( f_1, \ldots , f_{N_j} \right) ^T$ and the entries of the matrix $A^{(j)} \in  \mathbb{R}^{N_j \times N_j} $ are
\begin{equation*}
A^{(j)}_{ik}= \phi_{\varepsilon}^{(k)} (x_i,y_i), \quad i,k=1, \ldots, N_j.
\label{A}
\end{equation*}

If  the considered RBF is strictly positive definite, the interpolant \eqref{intg} is unique \cite{Fasshauer}. However, even if here we only focus on strictly positive definite RBFs, we remark that the uniqueness of the interpolant can be ensured also for the general case of strictly conditionally positive definite functions by modifying \eqref{rad1}, see e.g. \cite{Wendland05}.

In order to be able to formulate error bounds, we need some further assumptions on the regularity of $\Omega_j$ and thus we give the following definitions.
	
	\begin{definition}
		The fill distance, which is a measure of data distribution, is given by
		\begin{equation*}
		h_{ {\cal X}_N, \Omega} =  \sup_{ (x,y) \in \Omega} \left(  \min_{  (x_i,y_i)  \in {\cal X}_N} r_i(x,y) \right).
		\label{fd}
		\end{equation*}
	\end{definition}
	
		\begin{definition}
			An open and bounded covering $ \{ \Omega_j \}_{j=1}^{d}$ is called regular for $( \Omega, {\cal X}_N)$ if the following properties are satisfied \cite{Wendland02a}:
			\begin{itemize}
				\item[i.]  for each $ (x,y) \in \Omega$, the number of subdomains $ \Omega_j$ with $ (x,y) \in \Omega_j$ is bounded by a global constant $C$,
				\item[ii.]  there exist a constant $C_r > 0$ and an angle $\theta \in (0,\pi/2)$ such that every subdomain $ \Omega_j$ satisfies an interior cone condition (with angle $\theta$ and radius $ C_r h_{ {\cal X}_N, \Omega}$),
				\item[iii.]  the local fill distances $ h_{ {\cal X}_{N_j},\Omega_j}$ are uniformly bounded by the global fill distance $h_{{\cal X}_N, \Omega}$.
			\end{itemize}
			\label{constant}
		\end{definition}

		Thus, after defining the space $C_{ \nu}^{p}  ( \mathbb{R}^{2} ) $ of all functions $f \in C^p$ whose derivatives of order $ | \beta |=p $
		satisfy $ D^{ \beta} f ( x,y ) = {\cal O} ( (\sqrt{x^2+y^2})^{\nu} ) $ for $  \sqrt{x^2+y^2} \longrightarrow 0$, we consider the following convergence result \cite{Fasshauer,Wendland02a}.
		\begin{theorem}
			Let $ \Omega \subseteq  \mathbb{R}^2$ be open and bounded and suppose that $  {\cal X}_N= \{ (x_i,y_i)  ,i=1, \ldots, N \} \subseteq \Omega$. Let us consider a strictly conditionally  positive definite RBF which belongs to 
			$ C_{ \nu}^{p}  ( \mathbb{R}^{2} ) $. Let $ \{ \Omega_j \}_{j=1}^{d}$ be a regular covering for  $( \Omega,  {\cal X}_N)$  and let $ \{ W_j \}_{j=1}^{d}$ be $p$-stable for $ \{ \Omega_j  \}_{j=1}^{d}$. Then the error between $ f \in \mathscr{N}_{\phi} ( \Omega)$, where $ \mathscr{N}_{\phi}$ is the native space of the basis function,   and its PU interpolant \eqref{intg} can be bounded by
			$$
			| D^{ \beta} f( x,y ) -  D^{ \beta} {\cal I}( x,y) | \leq C^{'} h_{  {\cal X}_N, \Omega}^{\frac{ p+ \nu }{2} - | \beta |} |f|_{{\cal N}_{\phi} ( \Omega )},
			$$
			for all $(x,y) \in \Omega $ and all $ | \beta | \leq p/2$.
			\label{th1}
		\end{theorem}
		
\subsection{The positive constrained partition of unity approximant}

If we compare the result reported in Theorem \ref{th1} with the global error estimate shown in \cite{Fasshauer,Wendland05}, we can see that the PU interpolant preserves the local approximation order for the global fit. In particular, the PU method can be thought as the Shepard's method where $R_j$ are used instead of the data values $f_j$. Even if the classical Shepard's approximant in its original formulation is overcome \cite{Shepard}, it possesses a useful property; specifically, it lies within the range of the data. As a consequence, it is positive if the data values are positive \cite{Gordon}. The positivity-preserving property does not hold in its quadratic formulation nor for the PU approximant presented in the previous section.

In order to avoid such drawback for the PU method, we can directly act on the local RBF interpolants, following the strategy proposed in \cite{WU10a}. In such paper, a scheme devoted to construct a positive global fit is performed by defining several constraints. Extensive results show the good performances of such approach. In fact, the fit of positive samples is always positive, but, in comparison with the original unconstrained interpolation, a degrade of the quality of the approximation is observed. This is mainly due to the fact that a global method is considered and thus the shape of the surface is consequently globally modified (and not only in the area in which the interpolant is negative). Therefore, here we propose a new formulation for a positive fit considering \eqref{intg}.

Sufficient condition to have positive approximants on each subdomain $\Omega_j$ is that the coefficients $c_k$ of \eqref{rad1} are all positive. To such scope, following \cite{WU10a}, at first we choose ${\hat N}_j$  added data
\begin{equation*}
({\hat x}_{N_j+1},{\hat y}_{N_j+1}), \ldots, ({\hat x}_{N_j+{\hat N}_j},{\hat y}_{N_j+{\hat N}_j}),
\end{equation*} 
on the subdomain $\Omega_j$. Then, the $j$-th approximation problem consists in finding a function ${\hat R}_j$ of the form
\begin{equation} 
{\hat R}_j (x,y)= \sum_{k=1}^{N_j} c_k \phi_{\varepsilon}^{(k)} (x,y) + \sum_{{\hat k} = N_j+1}^{N_j+{\hat N}_j} c_{\hat k} \hat{\phi}_{\varepsilon_{{\hat k}}}^{({\hat k})}  (x,y), \quad (x,y) \in \Omega_j,
\label{rad_pos}
\end{equation}
such that
\begin{gather}
{\hat R}_j (x_i,y_i)= f_i,  \quad i=1, \ldots, N_j, \quad
c_i \geq 0, \quad i=1, \ldots , N_j +{\hat N}_j,
\label{con_pos}
\end{gather}
where $\hat{\phi}_{\varepsilon_{{\hat k}}}^{({\hat k})}$ are CSRBFs.

Note that in \eqref{rad_pos} we consider different supports for the CSRBFs. In particular, if a constraint $(\hat{x}_i,\hat{y}_i)$ is added  on $\Omega_j$ in a neighborhood of a point, namely $(x_i,y_i)$,  we select $\varepsilon_{{ \hat i}}$ such that only  $(x_i,y_i)$ belongs to the  support of the CSRBF. This choice is due to the fact that, doing in this way, at least when  $\hat{N}_j=N_j$, the problem \eqref{rad_pos} subject to \eqref{con_pos} admits solution. This is proved in \cite{WU10a}  by using the Gordan's Theorem \cite{Borwein}. A brief sketch of the proof will be given in what follows because important considerations arise. Let us define
\begin{gather*}
\boldsymbol{a}_k^{T}= - \{\phi_{\varepsilon}^{(k)}(x_i,y_i) \}_{i=1}^{N_j}, \quad k= 1, \ldots, N_j,  \nonumber\\
\boldsymbol{a}_{{\hat k}}^{T}= - \{ \hat{\phi}_{\varepsilon_{{\hat k}}}^{({\hat k})}(x_i,y_i) \}_{i=1}^{N_j} = (0, \ldots ,0, -b_{{\hat k}}, 0, \ldots, 0), \quad {\hat k}= N_j+1, \ldots ,2N_j,  \nonumber\\
\boldsymbol{a}_{2N_j+1}^{T}=(f_1, \ldots, f_{N_j}) , \nonumber
\label{2nj}
\end{gather*}
where $b_{{\hat k}}$ are positive real numbers. Then, since  a vector $\boldsymbol{v}$ such that $\boldsymbol{a}_i^{T} \boldsymbol{v}< 0$, $i=1, \ldots, 2N_j+1$, does not exist,  from Gordan's Theorem, we know that there exist nonnegative real numbers $c_1, \ldots ,c_{2N_j+1}$, such that $\sum_{i=1}^{2N_j+1} c_i \boldsymbol{a}_i=0$ and $c_{2N_j+1}>0$.

In practical use, the aim is to add as few data as possible,  therefore, according to \cite{WU10a}, one can find the minimum of $ \sum_{i=N_j+1}^{2N_j}
\textrm{sign}(c_{i})$, such that \eqref{con_pos} is satisfied.
However, this approach does not guarantee an optimal solution in terms of accuracy. Here instead, with a technique described in the next section, we select the optimal number of added data $\hat{N}_j$ which yields maximal accuracy. Nevertheless, our aim consists in perturbing \eqref{rad1} with \emph{small} quantities. A feasible way is to find 

	\begin{equation}
	\min \left( \sum_{i=N_j+1}^{N_j+\hat{N}_j} c_i^2 \right)^{1/2},
	\label{minprobl} 
	\end{equation}
	such that \eqref{con_pos} is satisfied.

The local approach here proposed enables us to modify the shape of the surface only if a negative local approximant is found, in fact if the $j$-th original local fit of the form \eqref{rad1} is positive, we do not need to add other data and in this case ${\hat R}_j=R_j$ on $\Omega_j$.  
Therefore, for each subdomain,  after selecting a  suitable number of constraints $\hat{N}_j$, which can also be 0,   the PC-PU approximant assumes the form
\begin{equation}
{\hat{\cal I}}(x,y)  = 
\sum_{j=1}^{d}  \left( \sum_{k=1}^{N_j} c_k \phi_{\varepsilon}^{(k)} (x,y) + \sum_{{\hat k} =  N_j+1}^{N_j+{\hat N}_j} c_{\hat k} \hat{\phi}_{\varepsilon_{{\hat k}}}^{({\hat k})} (x,y) \right) W_j (x,y). 
\label{intg_pos}
\end{equation}

\begin{remark}
	In \cite{WU10a} the authors limit their attention to CSRBFs. Here instead we will couple them with globally defined RBFs. In fact, even if specific supports of compactly supported kernels must be associated to the added constraints (see the second term in the right-hand side of \eqref{rad_pos}), we do not have any restrictions on the first term in the right-hand side of \eqref{rad_pos}. Thus, we can use different types of RBFs. In what follows, we will point out that coupling RBFs and CSRBFs leads to a benefit in terms of accuracy.
\end{remark}

\section{The PC-PU approximant: algorithm and suitable selection of positive constraints}
\label{locations}

This section is devoted to describe the PC-PU algorithm. In particular, we  focus on the suitable choice of the number of positive constraints.
Even if the technique here discussed is robust enough to work in any domain $\Omega \subseteq \mathbb{R}^2$, for simplicity we consider $\Omega=[0,1]^2$.

\subsection{Selection of positive constraints}

Acting as explained in the previous section, we can ensure the positivity of the PU approximant. However, depending on the number of positive constraints, this might lead to a low accuracy.  In \cite{WU10a}, for a global RBF-based interpolant, $h$ random data are selected as constraints if $f_h$ is the minimum among all values $f_i$, $i=1, \ldots, N$. This criterion does not guarantee a good accuracy neither the existence of a solution and thus we design an alternative approach  enabling us to select a suitable number of positive constraints.  

The proposed method is based on an a priori error estimate. To this aim, several schemes have already been developed. Precisely, we  focus on the so-called \emph{cross validation} algorithm, see \cite{Fasshauer,Golberg}. A variant of such method, known in literature as Leave One Out Cross Validation (LOOCV), is detailed in \cite{Rippa}. This approach is always used to find the optimal value of the shape parameter of the basis function. Here instead, for each PU subdomain we are interested in selecting a suitable number of constraints $ \hat{N}_j$. 

To avoid complexities, let us first consider an interpolation problem on $\Omega_j$ of the form \eqref{rad1}. Moreover, let us define the $j$-th interpolant obtained leaving out the $i$-th data on $\Omega_j$  as
\begin{equation*}
R^{i}_j (x,y)=  \sum_{k=1, k \neq i}^{N_j} c_k \phi_{\varepsilon}^{(k)} (x,y), \quad (x,y) \in \Omega_j,
\end{equation*}
and let
\begin{equation*}
e_i= f_i-R^{i}_j (x_i,y_i),
\end{equation*}
be the  error at the $i$-th point. Then, the quality of the local fit is determined by \emph{some} norm of the vector of errors $ ( e_1, \ldots ,e_{N_j})^{T}$, obtained by removing in turn one of the data points and comparing the resulting fit with the known value at the removed point. 

This implementation is computationally expensive. In fact, the matrix inverse, which requires ${\cal O} (N_j^3)$ operations, must be computed for each node. This leads to a total computational cost of ${\cal O} (N_j^4)$ operations. Thus, following \cite{Fasshauer,Rippa}, we simplify the computation to a single formula. Precisely, we calculate
\begin{equation}
e_i= \frac{c_i}{\left(A^{(j)}_{ii} \right)^{-1}},
\label{er0}
\end{equation}
where $c_i$ is the $i$-th coefficient of the interpolant based on the full data set and $(A^{(j)}_{ii})^{-1}$ is the $i$-th diagonal element of the inverse of the corresponding local interpolation matrix.

In our case, in order to guarantee a positive fit, we deal with an \emph{augmented} local problem and therefore, as error estimate, we compute the following quantity
\begin{equation}
\left( \hat{e}_1, \ldots, \hat{e}_{N_j + \hat{N}_j} \right)= \left( \frac{c_1}{\left(\hat{A}^{(j)}_{11}\right)^{-1}}, \ldots, \frac{c_{N_j + \hat{N}_j}}{\left(\hat{A}^{(j)}_{N_j + \hat{N}_j N_j+ \hat{N}_j}\right)^{-1}} \right),
\label{eq_man}
\end{equation}
where the symmetric matrix $\hat{A}^{(j)}$ is defined as
\begin{equation}
\hat{A}^{(j)}=
\begin{pmatrix}
\phi_{\varepsilon}^{(1)} (x_1,y_1) & \cdots & \hat{\phi}_{\varepsilon_{N_j+\hat{N}_j}}^{(N_j+\hat{N}_j)} ( x_{1},y_{1})\\
\vdots  & \ddots & \vdots \\
\hat{\phi}_{\varepsilon_{N_j+\hat{N}_j}}^{(N_j+\hat{N}_j)} ( x_{1},y_{1}) & \cdots & \hat{\phi}_{\varepsilon_{N_j+ \hat{N}_j}}^{(N_j+ \hat{N}_j)} ( \hat{x}_{N_j+ \hat{N}_j},\hat{y}_{N_j+ \hat{N}_j})\\
\end{pmatrix}. 
\label{mat}
\end{equation}
Moreover, in order to stress the dependence of the error  on $\hat{N}_j$, we use the notation
\begin{eqnarray}
\boldsymbol{\hat{e}}_j (\hat{N}_j) =\left( \hat{e}_1, \ldots, \hat{e}_{N_j + \hat{N}_j} \right).
\end{eqnarray}
Note that in our case the coefficients are not determined by directly computing the solution of a linear system, but they are found out by solving \eqref{minprobl}, subject to \eqref{con_pos}. Thus, to be more precise, we should refer to this criterion as LOOCV-like method, but in order to keep common notations, we will go on calling it simply LOOCV.
Indeed, we are able to fix a criterion which enables us to select a suitable number of positive constraints. Specifically, focusing on the maximum norm, we compute \eqref{eq_man} for $\hat{N}_j=1, \ldots, N_j$. 
Thus, if on $\Omega_j$ a negative fit is observed, we add $\hat{N}_j$ positive constraints if
\begin{equation}
||\hat{\boldsymbol{e}}_j(\hat{N}_j)||_{\infty} = \min_{k=1, \ldots,N_j} ||\hat{\boldsymbol{e}}_j(k)||_{\infty},
\label{er2}
\end{equation}
and the fit is positive. It is easy to see that we automatically ensure that the conditions \eqref{con_pos} are satisfied, indeed they are fulfilled at least for ${\hat N}_j=N_j$.

The ${\hat N}_j$ new added data can be placed randomly  within the $j$-th patch, but numerically we observe that selecting \emph{well distributed} points on $\Omega_j$ leads to a better accuracy. As a consequence, on the subdomain $\Omega_j$ of centre $(\bar{x}_j,\bar{y}_j)$ and radius $\delta$, we consider ${\hat N}_j$ positive constraints, distributed as the seeds on a sunflower head, i.e. defined as \cite{Swinbank,Vogel}
\begin{equation}
({\hat x}_{k}, {\hat y}_{k})=( \bar{x}_j + u_k \cos \eta_k, \bar{y}_j+ u_k \sin \eta_k),
\label{gp}
\end{equation}
$k=N_j+1, \ldots,  N_j+ {\hat N}_j$, where
\begin{equation*}
u_k= \delta \dfrac{ \sqrt{k- 1/2}}{\sqrt{{\hat N}_j-1/2}} \quad \textrm{and} \quad \eta_k=\dfrac{4 k \pi}{1+\sqrt5}.
\end{equation*}

\begin{figure}[ht!]
	\begin{center}
		\makebox[\textwidth]{
			\begin{minipage}[b]{0.5\textwidth}
				\centering
				\includegraphics[width=\textwidth]{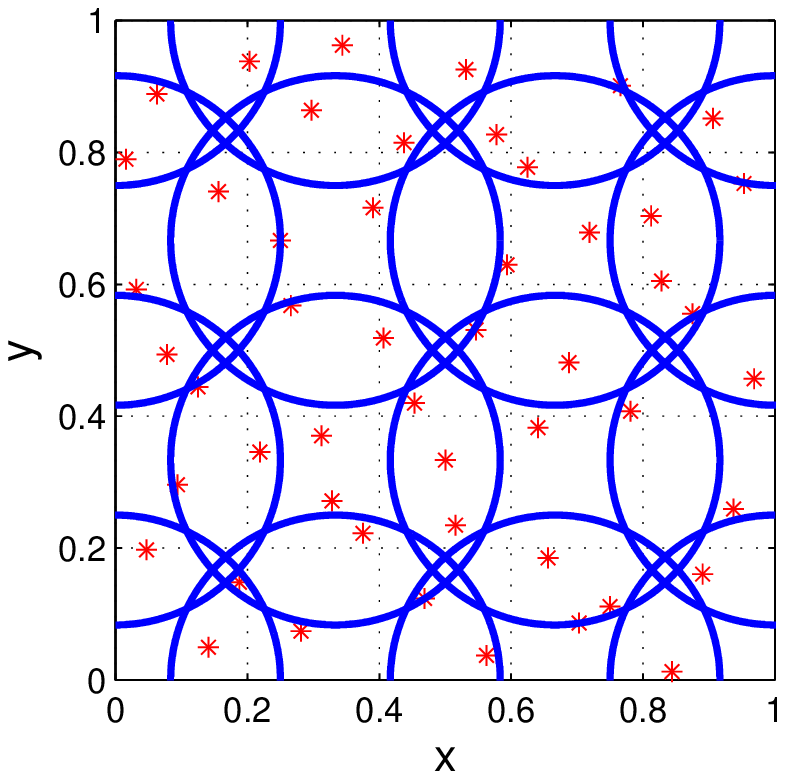} 
			\end{minipage}
			\begin{minipage}[b]{0.5\textwidth}
				\centering
				\includegraphics[width=\textwidth]{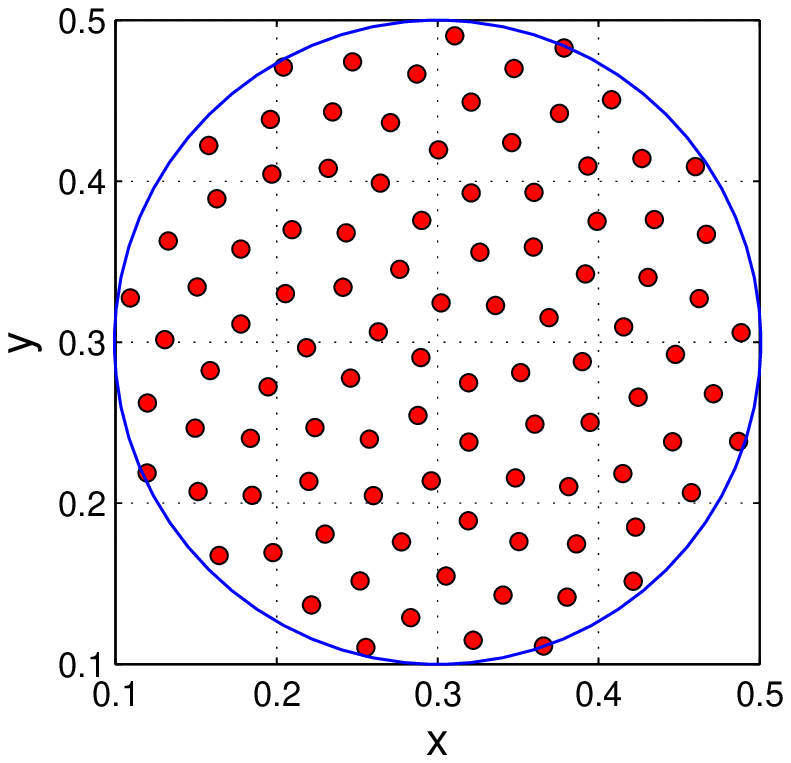} 
			\end{minipage}  
		}				
	\end{center}
	\caption{Left: an illustrative example of PU subdomains covering the domain $\Omega=[0,1]^2$. The red dots represent a set of scattered data and the blue circles identify the PU subdomains. Right: an illustrative example of added points computed with \eqref{gp}.}
	\label{figuraPUM1}
\end{figure}

We conclude this section with an illustrative figure.
A 2D view of a partition of unity structure covering a set  of scattered data in the unit square is shown in Figure \ref{figuraPUM1} (left); in the right frame we plot a set of points computed with \eqref{gp}. 

\subsection{Description of the algorithm}
To make simpler the presentation, we summarize in the {\tt PC-PU Algorithm} the steps needed to compute the PC-PU approximant.

At first, given a set of scattered data in $\Omega=[0,1]^2$,  we construct the PU subdomains. They are circular patches centred at a grid of points
$ {\cal C}_d = \{(\bar{x}_i, \bar{y}_i)   ,i=1, \ldots, d \} \subseteq \Omega$ of radius
\begin{equation}
\delta= \sqrt{\dfrac{1}{d}}.
\label{raggio2}
\end{equation}
According to \cite{Fasshauer}, we choose the number of  subdomains such that $N/d \approx 4$. The PC-PU approximant is then computed at a grid of $s$ evaluation points ${\cal E}_s= \{ (\tilde{x}_i,\tilde{y}_i)   ,i=1, \ldots, s \}$, see {\fontfamily{pcr} \selectfont Steps 1-3} of the {\tt PC-PU Algorithm}.

Once the partition of unity structure is built, a local interpolation problem needs to be solved for each PU subdomain. Specifically, in the $j$-th local approximation problem, only those data sites and evaluation points belonging to $\Omega_j$  are involved, see {\fontfamily{pcr} \selectfont Steps 4-6} of the {\tt PC-PU Algorithm}. To find such points the so-called kd-tree partitioning structures are commonly and widely used \cite{Arya98,Fasshauer}. However, here we use the efficient partitioning structure proposed in \cite{loocv_sottomess}. It is a multidimensional procedure which leads to a saving in terms of computational time with respect to the previous partitioning procedures proposed in \cite{Cavoretto15b,Cavoretto_CAMWA}.

Once both the sets $ {\cal X}_N$ and  ${\cal E}_s$  are organized among the different patches, an interpolant of the form \eqref{rad1} is constructed, see {\fontfamily{pcr} \selectfont Step 7} of the {\tt PC-PU Algorithm}. 
Then, only if the local fit is negative we add $\hat{N}_j$ positive constraints, as explained in the previous subsection, see {\fontfamily{pcr} \selectfont Steps 8-9} of the {\tt PC-PU Algorithm}. 
In this way, for each PU subdomain a positive local RBF approximant is computed for each evaluation point. Finally, the global fit \eqref{intg_pos} is applied accumulating all the $\hat{R}_j$ and $W_j$, see {\fontfamily{pcr} \selectfont Steps 10-11} of the {\tt PC-PU Algorithm}.

\begin{table}[ht!]
	\captionsetup{labelformat=empty}
	\begin{center}
		\begin{tabular}{p{12.5cm}*{1}{c}}
			\hline
			\vskip 0.01 cm 
			{\fontfamily{pcr} \selectfont INPUTS:} $N$, number of data; ${\cal X}_N=\{(x_i,y_i), i=1,\ldots,N\}$, set of data points;
			\vskip 0.08 cm 
			\hskip 2.2 cm   ${\cal F}_N=\{f_i, i=1,\ldots,N\}$, set of data values; $d$, number of 
			\vskip 0.08 cm 
			\hskip 2.2 cm subdomains;  $s$, number of evaluation points.
			\vskip 0.12 cm 
			{\fontfamily{pcr} \selectfont OUTPUTS:} ${\cal A}_s=\{ {\hat{\cal I}}(\tilde{x}_i,\tilde{y}_i), i=1,\ldots,s\}$, set of approximated values.
			\vskip 0.12 cm
			{\fontfamily{pcr} \selectfont Step 1:}
			A grid of  evaluation points ${\cal E}_s = \{(\tilde{x}_i,\tilde{y}_i), i=1,\ldots,s\} \subseteq \Omega$  is generated.
			\vskip 0.12 cm
			{\fontfamily{pcr} \selectfont Step 2:}
			A grid  of subdomain points ${\cal C}_d=\{(\bar{x}_j,\bar{y}_j), j=1,\ldots,d\} \subseteq \Omega$  is 
			\vskip 0.08 cm
			\hskip 2.2 cm  constructed.
			\vskip 0.12 cm
			{\fontfamily{pcr} \selectfont Step 3:} For each PU centre $(\bar{x}_j,\bar{y}_j)$, $j=1,\ldots,d$, a subdomain, whose
			\vskip 0.08 cm
			\hskip 2.2 cm  radius is given by \eqref{raggio2}, is constructed.
			\vskip 0.12 cm
			{\fontfamily{pcr} \selectfont Step 4:}  For each patch $\Omega_j$, $j=1,\ldots,d$,
			\vskip 0.12 cm
			\hskip 1.2 cm  {\fontfamily{pcr} \selectfont Step 5:} Find all data points ${\cal X}_{N_j}$  belonging to  $\Omega_j$.
			\vskip 0.12 cm
			\hskip 1.2 cm {\fontfamily{pcr} \selectfont Step 6:} Find all evaluation points ${\cal E}_{s_j}$ belonging to  $\Omega_j$.
			\vskip 0.12 cm
			\hskip 1.2 cm {\fontfamily{pcr} \selectfont Step 7:} Solve the unconstrained interpolation problem and a  local
			\vskip 0.08 cm
			\hskip 3.3 cm   interpolant $R_j$ is formed as in \eqref{rad1}.
			\vskip 0.12 cm
			\hskip 1.2 cm {\fontfamily{pcr} \selectfont Step 8:} If the local fit is positive $\hat{N}_j=0$ and $\hat{R}_j=R_j$,  else
			\vskip 0.12 cm
			\hskip 2.4 cm {\fontfamily{pcr} \selectfont Step 9:}
			For $\hat{N}_j=1 , \ldots, N_j$, compute the  constraints as in \eqref{gp} 
			\vskip 0.08 cm
			\hskip 4.4 cm  and calculate \eqref{eq_man}.
			\vskip 0.12 cm
			\hskip 3.3 cm  Consider ${\hat N}_j$ constraints, if ${\hat N}_j$ satisfies \eqref{er2} and \eqref{con_pos}.
			\vskip 0.12 cm
			\hskip 1.2 cm {\fontfamily{pcr} \selectfont Step 10:} Solve the constrained approximation problem and   form  a   
			\vskip 0.08 cm
			\hskip 3.3 cm  positive local approximant ${\hat R}_j$, see  \eqref{rad_pos}.
			\vskip 0.12 cm
			{\fontfamily{pcr} \selectfont Step 11:} 
			The interpolant \eqref{intg_pos} is formed by the weighted sum of the local fits.
			\\[\smallskipamount]
			\hline
		\end{tabular}
	\end{center}
	\caption{The {\tt PC-PU Algorithm}. Routine performing the PC-PU method.}
	\label{PUM_code_pos}
\end{table}

\section{Numerical experiments}
\label{NE_POS}
This section is devoted to show, by means of extensive numerical simulations,  the performances of the PC-PU approximant. 
In \cite{WU10a} the authors point out that the global CSRBF constrained method possesses a better approximation behaviour than the CMQS approximant \cite{Asim,Brodlie}.

Here, we compare our PC-PU fit with the original Shepard's method. In fact, for this method the positivity-preserving property holds, while it is lost in the modified quadratic version \cite{Gordon}. Then, comparisons with the global method proposed in \cite{WU10a} and with the classical PU interpolant will be carried out. Obviously, since we perform a PU approximation, large data sets are considered. On the opposite, a global interpolant, such as the one outlined in \cite{WU10a}, cannot handle large sets.

Experiments are performed considering several sets of random nodes contained in the unit square $\Omega = [0, 1]^2$, a grid of $d = \lfloor \sqrt{N} /2 \rfloor^2$ subdomain centres and a grid of $s = 80 \times 80$ evaluation  points.

In order to test the accuracy of the proposed method, we compute, for different kernels with different order of smoothness, the Maximum Absolute Error (MAE) and the Root Mean Square Error (RMSE) whose formulas are
\begin{align*} 
MAE = \max_{1\leq i \leq s} |f(\tilde{x}_i,\tilde{y}_i) - {\cal \hat{I}}(\tilde{x}_i, \tilde{y}_i)|, 
\end{align*}
and
\begin{align*} 
RMSE = \sqrt{\frac{1}{s}\sum_{i=1}^{s} |f(\tilde{x}_i,\tilde{y}_i) - {\cal \hat{I}}(\tilde{x}_i, \tilde{y}_i)|^2}.
\end{align*}
The errors are  computed using  as test functions
\begin{align*}
f_1(x,y)=(x-0.5)^2+(y-0.4)^2,       
\end{align*}
and
\begin{align*}
f_2(x,y)=   [3 (y-0.4) \sin (x-0.5)]^2(y+0.5)^{1/3}.        
\end{align*}

Experiments  are carried out considering  the Wendland's $C^2$ function for the added nodes in \eqref{rad_pos}. Furthermore, since for the given interpolation conditions the choice is arbitrary, we use both the  Wendland's $C^2$ and the IMQ $C^{\infty}$ functions.  For them, we have to fix the shape parameters. We remark that the results are affected by the choice of the $\varepsilon$, in fact small values of $\varepsilon$ lead to problems of instability, while for large values the approximate solution might be inaccurate. In particular, many researchers already worked on the problem of finding stable approximations  when $\varepsilon$ tends to zero (see e.g. \cite{Demarchi15,Driscoll-Fornberg02,Fornberg11}). Thus, referring to such papers one can easily guess which are the \emph{safe} values. Here, for instance,  we fix the shape parameter $\varepsilon$ equal to $0.1$ and $1$ for the Wendland's $C^2$ and the IMQ $C^{\infty}$ functions, respectively. Moreover, we point out that our aim is to compare the behaviour of the PC-PU method versus the classical PU approach. In this sense, such comparison is independent from the shape parameter, indeed we always register a decrease of the accuracy when the PC-PU method is used.

\begin{figure}[h]
	\begin{center}
		\makebox[\textwidth]{
			\begin{minipage}[b]{0.5\textwidth}
				\centering
				\includegraphics[width=\textwidth]{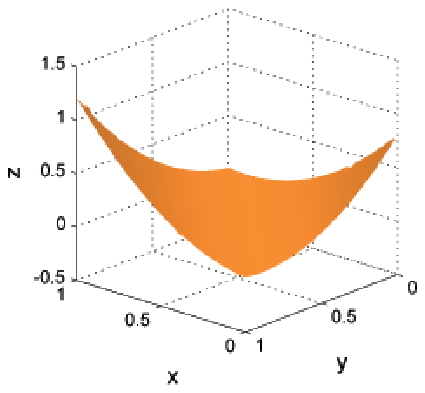} 
			\end{minipage}
			\begin{minipage}[b]{0.5\textwidth}
				\centering
				\includegraphics[width=\textwidth]{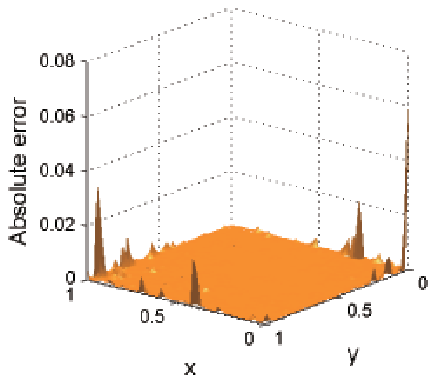} 
			\end{minipage}  
		}				
	\end{center}
	\begin{center}
		\makebox[\textwidth]{
			\begin{minipage}[b]{0.5\textwidth}
				\centering
				\includegraphics[width=\textwidth]{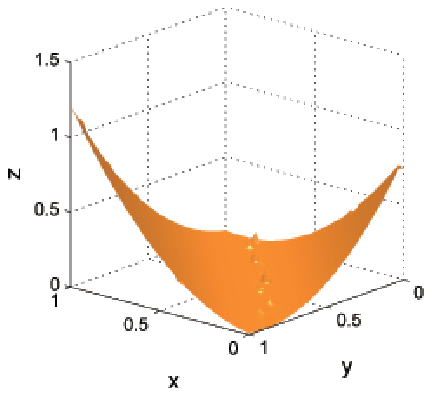} 
			\end{minipage}
			\begin{minipage}[b]{0.5\textwidth}
				\centering
				\includegraphics[width=\textwidth]{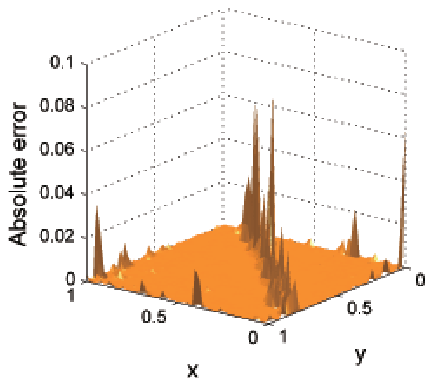} 
			\end{minipage}  
		}				
	\end{center}			
	\caption{The function $f_1$ (left) and absolute errors (right) obtained by applying the PU (top) and PC-PU (bottom) with the  Wendland's $C^2$ kernel.} 
	\label{fig_PCPUCS}
\end{figure}

\begin{figure}[h]
	\begin{center}
		\begin{center}
			\makebox[\textwidth]{
				\begin{minipage}[b]{0.5\textwidth}
					\centering
					\includegraphics[width=\textwidth]{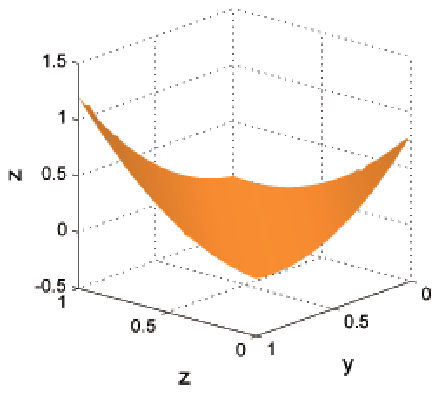} 
				\end{minipage}
				\begin{minipage}[b]{0.5\textwidth}
					\centering
					\includegraphics[width=\textwidth]{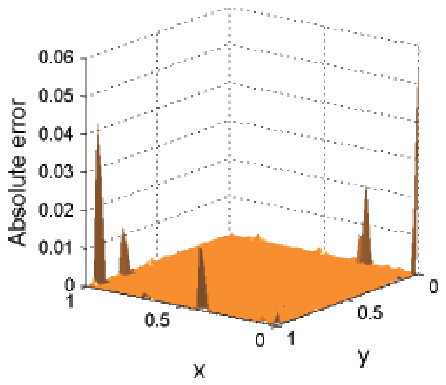} 
				\end{minipage}  
			}				
		\end{center}
		\begin{center}
			\makebox[\textwidth]{
				\begin{minipage}[b]{0.5\textwidth}
					\centering
					\includegraphics[width=\textwidth]{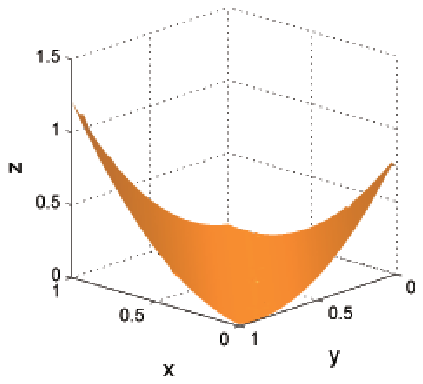} 
				\end{minipage}
				\begin{minipage}[b]{0.5\textwidth}
					\centering
					\includegraphics[width=\textwidth]{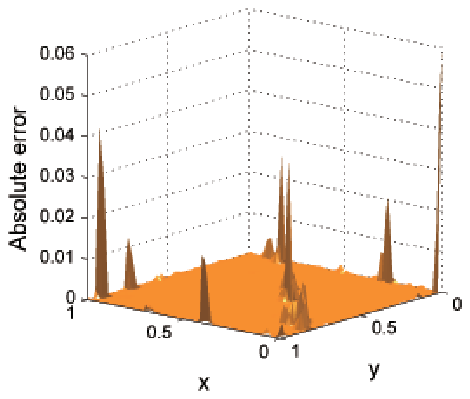} 
				\end{minipage}  
			}				
		\end{center}		
		\caption{The function $f_1$ (left) and absolute errors (right) obtained by applying the PU (top) and PC-PU (bottom) with the IMQ $C^\infty$ kernel.} 
		\label{fig_PCPU}
	\end{center}
\end{figure}

Tables \ref{ta1} and \ref{ta2} show a direct comparison between the classical PU interpolant, which leads to a negative fit,  and the PC-PU approximant.
We can observe a better behaviour of the PC-PU approximant when globally supported RBFs are used instead of CSRBFs. In fact, with the latter, large errors are introduced in the region where a negative fit is observed. Roughly speaking, using globally defined RBFs, the errors of the classical PU method and of the PC-PU approach are close to each other. On the opposite, by means of CSRBFs, the PC-PU approximant preserves the positivity property with an error which is about two times the one of the unconstrained interpolant. In order to have a graphical prove, refer to Figures \ref{fig_PCPUCS} and \ref{fig_PCPU} in which we consider $N=3500$ nodes and the test function $f_1$. From these figures we can note that, both with the use of CSRBFs and globally defined RBFs there is no smoothing effect when the solution approaches zero, i.e. when constraints are used. 

\begin{figure}[ht!]
	\begin{center}
		\makebox[\textwidth]{
			\begin{minipage}[b]{0.5\textwidth}
				\centering
				\includegraphics[width=\textwidth]{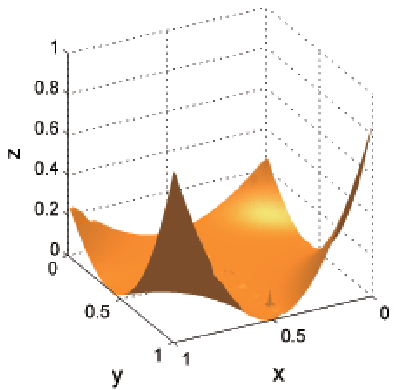} 
			\end{minipage}
			\begin{minipage}[b]{0.5\textwidth}
				\centering
				\includegraphics[width=\textwidth]{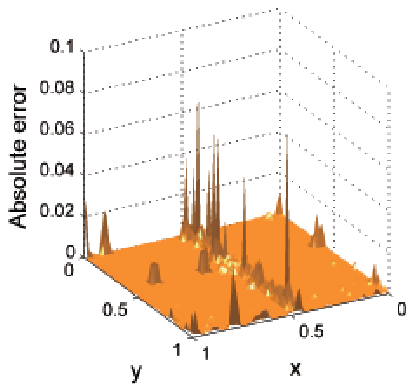} 
			\end{minipage}  
		}				
	\end{center}
	\begin{center}
		\makebox[\textwidth]{
			\begin{minipage}[b]{0.5\textwidth}
				\centering
				\includegraphics[width=\textwidth]{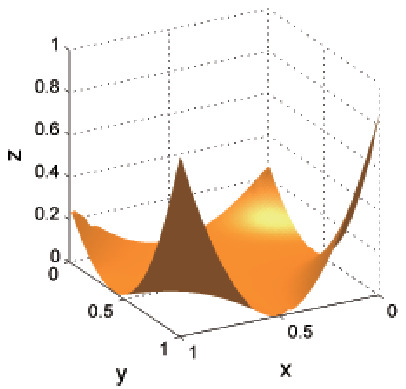} 
			\end{minipage}
			\begin{minipage}[b]{0.5\textwidth}
				\centering
				\includegraphics[width=\textwidth]{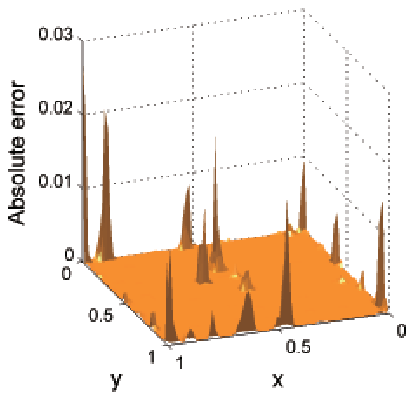} 
			\end{minipage}  
		}				
	\end{center}	
	\caption{The function $f_2$ (left) and absolute errors (right) obtained by applying  the  PC-PU method with the Wendland's $C^2$ (top) and  the IMQ $C^\infty$ (bottom) kernels}. 
	\label{fig_con}
\end{figure}

\begin{figure}[ht!]
	\begin{center}
		\makebox[\textwidth]{
			\begin{minipage}[b]{0.5\textwidth}
				\centering
				\includegraphics[width=\textwidth]{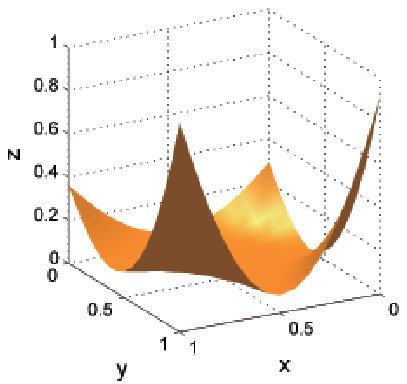} 
			\end{minipage}
			\begin{minipage}[b]{0.5\textwidth}
				\centering
				\includegraphics[width=\textwidth]{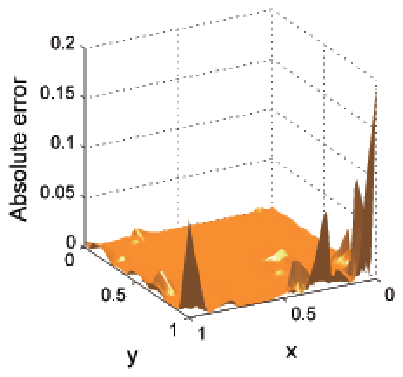} 
			\end{minipage}  
		}				
	\end{center}
	\begin{center}
		\makebox[\textwidth]{
			\begin{minipage}[b]{0.5\textwidth}
				\centering
				\includegraphics[width=\textwidth]{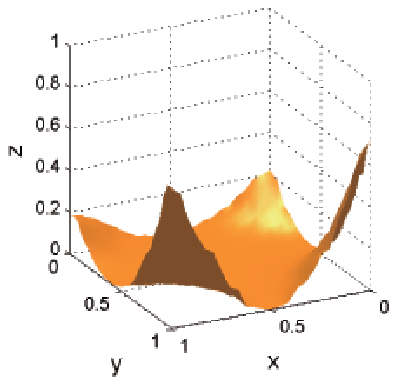} 
			\end{minipage}
			\begin{minipage}[b]{0.5\textwidth}
				\centering
				\includegraphics[width=\textwidth]{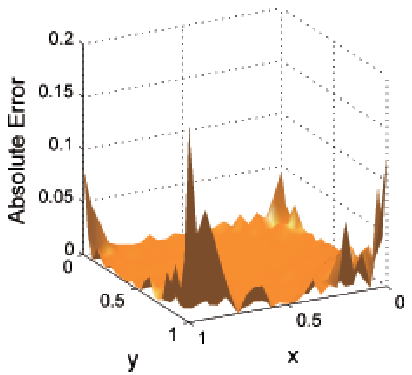} 
			\end{minipage}  
		}				
	\end{center}			
	\caption{The function $f_2$ (left) and absolute errors (right) obtained by applying   the CSRBF-based global method (top) and the Shepard's approximant (bottom).} 
	\label{fig_con1}
\end{figure}

Then, in Tables \ref{ta3} and \ref{ta4} we show the errors obtained with the Shepard's method and the CSRBF-based global method  proposed in \cite{WU10a}. For this method we consider the Wendland's $C^2$ kernel as CSRBF. 
As evident from Figures \ref{fig_con} and \ref{fig_con1}, in which we consider $N=3500$ random nodes and the test function $f_2$, our proposed local scheme maintains a better accuracy than the other considered methods. Moreover, from Figure \ref{fig_con}, we can again observe a better behaviour of the PC-PU method coupled with globally defined RBFs rather than with CSRBFs. This behaviour, is also due to the fact that with CSRBFs we need to add more constraints. In particular, if on $\Omega_j$ a negative fit is observed, the PC-PU technique selects ${\hat N}_j$ constraints, computed with \eqref{gp}, such that the  positivity of the local interpolant is ensured. In general, we note that when globally defined RBFs are used the number ${\hat N}_j$ is so that ${\hat N}_j \ll N_j$. On the opposite with CSRBFs, ${\hat N}_j$ is usually closer to $N_j$. In other words, we need to add more constraints in case of CSRBFs and, as stressed in \cite{WU10a}, this causes a decrease of the fit accuracy.

\begin{table}[h]
	\begin{center}
		\begin{tabular}{cccc} 	\hline\noalign{\smallskip}
			$N$ & method	 & MAE & RMSE  \\
			\noalign{\smallskip}
			\hline
			\noalign{\smallskip}
			$\hskip-2pt 300$ & PU &  $1.50{\rm E}-01$ & $1.52{\rm E}-02$    	   \\	
			&  PC-PU      &  $1.50{\rm E}-01$   & $2.03{\rm E}-02$       	   \\
			\rule[0mm]{0mm}{3ex}		
			$1000$     &  PU   &  $7.36{\rm E}-02$   & $3.08{\rm E}-03$     \\
			&   PC-PU  &  $7.96{\rm E}-02$ & $6.44{\rm E}-03$      \\
			\rule[0mm]{0mm}{3ex}
			$3500$     & PU  &  $6.34{\rm E}-02$ & $1.47{\rm E}-03$    	   \\
			&   PC-PU    &  $8.40{\rm E}-02$ 	  & $2.86{\rm E}-03$        \\
			\rule[0mm]{0mm}{3ex}
			$8000$     &  PU  &  $2.43{\rm E}-02$    & $4.17{\rm E}-04$    \\
			&   PC-PU   &  $5.99{\rm E}-02$   & $1.03{\rm E}-03$         \\
			\hline 
		\end{tabular}
	\end{center}
	\caption{MAEs and RMSEs  computed with the Wendland's $C^2$ kernel for $f_1$.}
	\label{ta1}
\end{table}

\begin{table}[h]
	\begin{center}
		\begin{tabular}{cccc} 	\hline\noalign{\smallskip}
			$N$ & method	 & MAE & RMSE  \\
			\noalign{\smallskip}
			\hline
			\noalign{\smallskip}
			$\hskip-2pt 300$ & PU &  $1.39{\rm E}-01$ & $1.04{\rm E}-02$    	   \\	
			&  PC-PU      &  $1.39{\rm E}-01$   & $1.44{\rm E}-02$       	   \\
			\rule[0mm]{0mm}{3ex}
			$1000$     &  PU   &  $7.02{\rm E}-02$   & $2.88{\rm E}-03$     \\
			&   PC-PU  &  $7.02{\rm E}-02$ & $3.49{\rm E}-03$      \\
			\rule[0mm]{0mm}{3ex}
			$3500$     & PU  &  $5.89{\rm E}-02$ & $1.50{\rm E}-03$    	   \\
			&   PC-PU    &  $5.89{\rm E}-02$ 	  & $1.66{\rm E}-03$        \\
			\rule[0mm]{0mm}{3ex}
			$8000$     &  PU  &  $2.33{\rm E}-02$    & $3.46{\rm E}-04$    \\
			&   PC-PU   &  $2.33{\rm E}-02$   & $3.68{\rm E}-04$         \\
			\hline 
		\end{tabular}
	\end{center}
	\caption{MAEs and RMSEs  computed with the IMQ $C^{\infty}$ kernel for $f_1$.}
	\label{ta2}
\end{table}

\begin{table}[h]
	\begin{center}
		\begin{tabular}{cccc} 	\hline\noalign{\smallskip}
			$N$ & method	 & MAE & RMSE  \\
			\noalign{\smallskip}
			\hline
			\noalign{\smallskip}
			$\hskip-2pt 300$ & CSRBF global &  $3.36{\rm E}-01$ & $2.77{\rm E}-01$    	   \\	
			&  PC-PU      &  $2.76{\rm E}-01$   & $1.91{\rm E}-02$       	   \\
			\rule[0mm]{0mm}{3ex}
			$1000$     &   CSRBF global   &  $4.31{\rm E}-01$   & $2.59{\rm E}-02$     \\
			&   PC-PU  &  $8.84{\rm E}-02$ & $5.95{\rm E}-03$      \\
			\rule[0mm]{0mm}{3ex}
			$3500$     &  CSRBF global  &  $1.69{\rm E}-01$ & $6.67{\rm E}-02$    	   \\
			&   PC-PU    &  $8.48{\rm E}-02$ 	  & $2.61{\rm E}-03$        \\
			\hline 
		\end{tabular}
	\end{center}
	\caption{MAEs and RMSEs of the CSRBF-based global method and of the PC-PU approximant computed with the Wendland's $C^2$ kernel for $f_2$.}
	\label{ta3}
\end{table}

\begin{table}[h]
	\begin{center}
		\begin{tabular}{cccc} 	\hline\noalign{\smallskip}
			$N$ & method	 & MAE & RMSE  \\
			\noalign{\smallskip}
			\hline
			\noalign{\smallskip}
			$\hskip-2pt 300$ & Shepard &  $8.38{\rm E}-01$ & $3.21{\rm E}-01$    	   \\	
			&  PC-PU      &  $1.32{\rm E}-01$   & $1.48{\rm E}-02$       	   \\
			\rule[0mm]{0mm}{3ex}
			$1000$     &   Shepard  &  $5.71{\rm E}-01$   & $7.75{\rm E}-02$     \\
			&   PC-PU  &  $8.62{\rm E}-02$ & $4.31{\rm E}-03$      \\
			\rule[0mm]{0mm}{3ex}
			$3500$     &  Shepard  &  $1.87{\rm E}-01$ & $1.42{\rm E}-02$    	   \\
			&   PC-PU    &  $2.89{\rm E}-02$ 	  & $9.73{\rm E}-04$        \\
			\hline 
		\end{tabular}
	\end{center}
	\caption{MAEs and RMSEs of the Shepard's method and of the PC-PU approximant computed with the IMQ $C^{\infty}$ kernel for $f_2$.}
	\label{ta4}
\end{table}

\section{Application to population dynamics}
\label{app_erbivori}

In this section we point out how the PC-PU procedure can be useful in applied sciences. Among several applications which can be investigated, we focus on a mathematical model dealing with wild herbivores in forests.

Let $H$, $G$ and $T$ represent respectively the herbivores, grass and trees populations of the environment in consideration. The model we consider  is a classical predator with two prey system \cite{TambVent}, in which the resources are consumed following a concave response function, usually called
the Beddington-De Angelis function \cite{De Angelis}.  Let $\alpha$ and $\beta$ be the inverse of the herbivores maximal consumption of grass  and trees, respectively. Let $r_1$ and $r_2$ denote the grass and trees growth rates, $K_1$ and $K_2$ their respective carrying capacities, $\mu$ the metabolic rate of herbivores, $c$ and $g$ the half saturation constants, $e\leq 1$ and $f\leq 1$ the conversion factors of food into new herbivore biomass and $a$ and $b$ the daily feeding rates due to grass and trees, respectively, the model  reads as follows \cite{Sabetta}
\begin{equation*}
\begin{array}{ll}
\frac{ \displaystyle  dH}{ \displaystyle  dt}=-\mu H +ae\dfrac{H G}{c+H+\alpha G}+bf\dfrac{H T}{g+H+\beta\ T+\alpha G},  & \textrm{} \\
\vspace{.01cm}\\
\frac{ \displaystyle  dG}{ \displaystyle  dt}=r_1 G\left( 1-\dfrac{G}{K_1}\right)-a\dfrac{H G}{c+H+\alpha G}, & \textrm{} \\
\vspace{.01cm}\\
\frac{ \displaystyle  dT}{ \displaystyle  dt}=r_2 T\left( 1-\dfrac{T}{K_2}\right)-b\dfrac{H T}{g+H+\beta T+\alpha G}.  & \textrm{} 
\end{array}
\label{HGTsys}
\end{equation*}
All parameters are nonnegative. In particular, $K_1$, $K_2$, $c$ and $g$ are measured in biomass,  $e$, $f$, $\alpha$ and $\beta$ are pure numbers, $\mu$, $r_1$, $r_2$, $a$ and $b$ are rates.

For the study of critical points and their stability, refer to \cite{TambVent}. The equilibrium that play a role in this investigation is  the coexistence equilibrium point $E^*=(H^*,G^*,T^*)$. It can be assessed only via numerical simulations. We test the model considering the natural park of the Dolomiti Bellunesi located in the  Northern Italy.  To estimate the parameters, we refer to \cite{Fuji,TambVent}.  Basing our considerations on data related to the number of herbivores, the extension of the park and following  tables providing the estimation of the annual net primary production of several environments  shown in \cite{Fuji}, we can fix $\mu=0.03$,  $r_1=0.01$, $r_2=0.0006$, $\alpha=0.05^{-1}$, $\beta=8$, $e=0.605$, $f=0.001$, $K_1={3469640.64}$, $K_2={15695993.39}$, $c={101862.16}$ and $g={1001229580.18}$. Moreover, also the initial conditions $H(0)=268.750$, $G(0)=2313093.76$ and $T(0)=1046399.56$ are known. 

\begin{figure}[ht!]
	\begin{center}
		\makebox[\textwidth]{
			\includegraphics[height=.24\textheight]{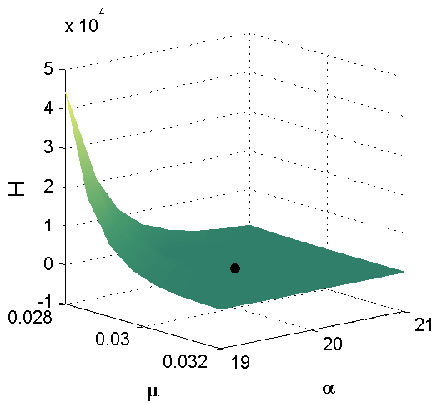} 
			\includegraphics[height=.24\textheight]{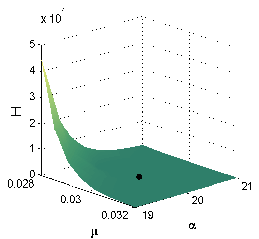}}
		\caption{Dolomiti Bellunesi's system. The surface represents the equilibrium value of the herbivores as function of
			the parameters $\mu$ and $\alpha$. The surface is reconstructed with the PU and PC-PU approximant, left and right respectively.}
		\label{applicazione_sabetta}
	\end{center}
\end{figure}

By means of numerical simulations carried out in \cite{Sabetta}, we can observe that the herbivore population level appears to be very sensitive to small perturbations of several parameters and under the  high risk of extinction. The results indicate that the parameters most affecting the system's final configuration are $\alpha$ and $\mu$.  Such consideration follows from the surface plotted in Figure \ref{applicazione_sabetta}; the surface shows the value attained by the herbivores in function of the parameters $\alpha$ and $\mu$.  The dot represents the situation in the present ecosystem conditions. As we can graphically note from the surface, for small perturbations on the current values of  $\alpha$ and $\mu$, the herbivores can extinguish, i.e. $H=0$.

In Figure \ref{applicazione_sabetta} the surfaces are reconstructed with the standard PU interpolant and with the PC-PU approximant, left and right respectively.
Even if it is not evident because of the large scale of the $z$-axis, the surface reconstructed  with the PU interpolant reaches its minimum at $-1$, violating the biological constraint $H \geq 0$. From this consideration, the importance of having a robust tool enabling us to fit positive data values with a positive approximant follows.

\section{Concluding remarks}
\label{cr}
In this paper we presented a robust technique devoted to preserve the positivity of the fit via a local RBF-based method. Extensive numerical simulations have been carried out to show the effectiveness of the method. 
Moreover, in order to motivate the importance of a meshless method that do not violate biological or physical constraints, we investigated an application in biomathematics.

Work in progress consists in extending the proposed tool in higher dimensions. Furthermore, in case of ill-conditioning, since in this context CSRBFs are less performing than globally defined RBFs, further investigations in the use of stable bases	\cite{Demarchi15,Fornberg11} are needed.

\end{document}